\documentclass[amsmath,twocolumn,secnumarabic,superscriptaddress,%
    floatfix,amssymb,aps]{revtex4}
\usepackage{times,graphicx,amsthm,url}
\usepackage[percent]{overpic}

\graphicspath{{./figs/}}
\newcommand{\K}{\mathcal{K}}
\newcommand{\bd}{\partial}
\newcommand{\wall}{\mathrm{Wall}}

\newcommand{\VF}{\mathrm{VF}}
\newcommand{\gStrut}{\mathrm{GStrut}}
\newcommand{\gThi}{\mathrm{GThi}} %Gehring thickness.
\newcommand{\len}{\mathrm{Len}}
\newcommand{\eps}{\varepsilon}
\newcommand{\dplus}{\delta^+}
\renewcommand{\d}{\mathrm{d}}
\newcommand{\disj}{\sqcup} %not coprod, which is an upside down Pi
\newcommand{\bigdisj}{\bigsqcup}
\newcommand{\grad}{\nabla}
\newcommand{\R}{\mathbb{R}}
\newcommand{\lem}[1]{Lemma~\ref{lem:#1}}
\newcommand{\figr}[1]{Figure~\ref{fig:#1}}
\newcommand{\thm}[1]{Theorem~\ref{thm:#1}}

\begin{document}

%\newcommand{\qed}{\par\vskip-\baselineskip\prbox\par}
%\setcounter{chapter}{0}
%\chapter{SOME ROPELENGTH-CRITICAL CLASPS}
%\markboth{Sullivan and Wrinkle}{Some Ropelength-Critical Clasps}
%\author{John~M.~Sullivan and Nancy~C.~Wrinkle}

%\address{Institut f\"ur Mathematik, Technische Universit\"at Berlin, MA 3-2,\\
%Str.~des 17.~Juni 136, 10623 Berlin, Germany \\
%and \\
%Department of Mathematics, University of Illinois,\\
%1409 W.~Green St., Urbana IL 61801, USA\\
%E-mail: jms@isama.org\\ \smallskip
%Department of Mathematics, Northeastern Illinois University,\\
%Science Building 214C, 5500 N.~St.~Louis Av., Chicago IL 60625, USA\\
%E-mail: n-wrinkle@neiu.edu}

%%%
\title{Some Ropelength-Critical Clasps}
\date{March 15, 2004; revised September 20, 2004}

\author{John M. Sullivan}
\email[Email: ]{jms@isama.org}
\affiliation{Institut f\"ur Mathematik, Technische Universit\"at Berlin,
D--10623 Berlin}
\affiliation{Department of Mathematics, University of Illinois,
Urbana, IL 61801}

\author{Nancy C. Wrinkle}
\email[Email: ]{n-wrinkle@neiu.edu}
\affiliation{Department of Mathematics, Northeastern Illinois University,
Chicago, IL 60625}

\newtheorem{theorem}{Theorem}[section]
\newtheorem{lemma}[theorem]{Lemma}
\newtheorem{corollary}[theorem]{Corollary}
\newtheorem{proposition}[theorem]{Proposition}
\theoremstyle{definition}
\newtheorem*{definition}{Definition}
%%%

\begin{abstract}
We describe several configurations of clasped ropes
which are balanced and thus critical for the Gehring ropelength problem.
\end{abstract}
%%%
\maketitle
%%%

\section{Introduction}     %S1-Heads

The ropelength of a link is given by the ratio of its length to its thickness.
There are many ways to measure the thickness of a space curve, but for links
one particularly simple notion, the \emph{Gehring thickness}, is simply the
minimum distance between different components.  With Cantarella, Fu and Kusner
we introduced a theory of criticality~[2] for the Gehring ropelength
problem.  Our necessary and sufficient conditions for ropelength criticality
take the form of a balance criterion which says that the tension force trying
to reduce the length of the curves must be balanced by contact forces acting
at points achieving the minimum distance.  One simple example, with surprising
intricacy for its solution~[2], is the \emph{clasp}.  This is
a generalized link whose components are not closed curves but instead have
constrained endpoints.  In the clasp, one rope whose ends are attached to
the ceiling is looped around another whose ends are attached to the floor.
The ropelength-critical clasp we described is presumably
the minimizer for ropelength, and is surprising in several ways:
the tips of the two components are 6\% further apart than they need to
be---leaving a small gap between the ropes---and the curvature
of the core curves blows up at the tips.

Here, we describe critical configurations of several generalized clasps with
two or more components.  Our new examples include clasps with one or both
curves doubled, and connect sums of clasps with Hopf links.
All of our examples are proven critical by the balance criterion.
Many of them we expect are global minimizers (though we know no way
to show this), but others are clearly not minimizing and are presumably
unstable equilibria.

\section{Background}
Here we recall the necessary definitions and theorems from
our work~[2] with Cantarella, Fu and Kusner.
These will show the balance criterion for Gehring ropelength
in the form in which we apply it to our new examples.

\begin{definition}
A \emph{generalized link} $L$ is a curve $L$ (with disjoint components)
together with obstacles and endpoint constraints.
In particular, each endpoint $x\in\partial L$ is constrained to stay on
some affine subspace $M_x\subset\R^3$.
Furthermore, there is a finite collection of \emph{obstacles} for the link,
each obstacle $$\{p\in\R^3: g_j(p)<0\}$$ being given in terms of a
$C^1$ function $g_j$ having $0$ as a regular value.
By calling these sets obstacles, we mean that $L$ is constrained
to stay in the region where $\min g_j \ge 0$.
\end{definition}

\begin{definition}
The \emph{Gehring thickness} $\gThi(L)$ of a curve $L$ is the minimal distance
between points on different components of $L$.
This is the supremal~$\tau$ for which the $(\tau /2)$-neighborhoods
of the components of $L$ are disjoint.
\end{definition}

We formulate the ropelength problem as minimizing the length of
a generalized link subject to the constraint that its Gehring
thickness remains at least~$1$.
The contact points of the different components of the link
are of primary importance, and are called \emph{struts}.

\begin{definition}
\label{defn:gstrut}
An (unordered) pair of points $x$ and $y$ on different components of $L$ is a
\emph{Gehring strut} if $|y-x| = \gThi(L)$.
The set of all Gehring struts of $L$ is denoted $\gStrut(L)$.
\end{definition}

Given a generalized link $L$, only variations preserving
the endpoint constraints should be allowed.
A continuous vectorfield $\xi$ along $L$ is said to be \emph{compatible}
with these constraints if it is tangent to $M_x$ at each endpoint $x\in\bd L$.
We write $\VF_c(L)$ for the space of all compatible vectorfields.

Given a set of obstacles $g_j<0$ for a link $L$, we write
$$O(L) := \min_j \min_{x\in L} g_{j(x)}.$$
Then $L$ avoids the obstacles $g_j$ if and only if $O(L)\ge 0$.
We define the wall struts of $L$ by
$$\wall_j(L) := L\cap\{g_j=0\},\quad \wall(L)=\bigdisj_j \wall_j(L).$$
This incorporates those parts of $L$ which are on the boundary of the obstacles.

Usually we will be minimizing the length~$\len(L)$ of a link while
constraining $\gThi(L)$ to be at least~$1$.  Sometimes, however, we
wish to consider a slightly more general objective functional,
a weighted sum $\len^w(L) := \sum w_i \len(L_i)$
of the lengths of the different components.  Here $w_i\ge 0$ can
be viewed as the elastic tension within component~$L_i$.

It is known~[2] that critical links for this (weighted) Gehring
ropelength problem are curves of finite total curvature.
(See~[3] for an expository account of such curves.)
This means that the first variation of length under a compatible $\xi$ is
given by $\delta_\xi \len(L) = -\int_L \left< \xi, \K \right>$, where~$\K$
is a vector-valued Radon measure along~$L$ (what we call a \emph{force}
along $L$) called the curvature force.  For a~$C^2$ link,
we have $\K=\kappa N\,ds$ in terms of the Frenet frame.
For weighted length, it follows that $\delta \len^w = -\K^w$ where
$\K^w=w_i\K$ along component~$L_i$.

A variation of~$L$ will change the length of the struts
and change the values of the obstacle functions~$g_j$.  We collect
these changes into the \emph{rigidity operator}~$A$.
If $L$ is varied with initial velocity $\xi\in\VF_c(L)$,
then $A(\xi)$ is by definition a continuous function on
$\gStrut(L)\disj\wall(L)$; its value on a strut $\{x,y\}$
of length $1$ is $\delta_\xi|x-y|=(\xi_x-\xi_y)\cdot(x-y)$,
and its value on a wall strut~$x$ (where $g_j(x)=0$) is
$\delta_\xi g_j(x)=\xi_x\cdot\grad g_j$.

The (one-sided) first variation $\dplus_{\xi} \gThi(L)$
of Gehring thickness is (by Clarke's differentiation
theorem for min-functions, see~[2]) the minimum of $A(\xi)$
over all struts; similarly the first variation $\dplus_{\xi} O(L)$
of the obstacle function is the minimum over all wall struts.

Intuitively, we expect a ropelength-critical configuration
to be one whose length cannot be reduced without also
reducing thickness.  For technical reasons (see~[2])
we define strong criticality to require this reduction
to happen at a definite rate:

\begin{definition}
\label{defn:strongcrit}
We say that a generalized link $L$ is
\emph{strongly critical} for minimizing weighted length
when constrained by $\gThi$
if there is an $\eps>0$ such that, for all compatible $\xi$
with $\delta_\xi \len^w = -1$, we have
$$\min\big(\dplus_{\xi} \gThi(L), \,\dplus_{\xi} O(L) \big)\le -\eps.$$
\end{definition}

The adjoint $A^*$ of the rigidity operator takes a Radon measure
on $\gStrut\disj\wall$ and gives a force along $L$.
This adjoint is what appears in our balance
criterion~[2], which in turn allows us to explicitly solve
for the shapes of critical configurations of various links.

\begin{theorem}\label{thm:gen-gehr-bal}
A generalized link $L$ is strongly critical for weighted
Gehring ropelength
if and only if there is a positive Radon measure $\mu$ on
$\gStrut(L) \disj \wall(L)$ such that $-\K^w=A^*\mu$
as linear functionals on $\VF_c(L)$.
That is, $-\K^w$ and $A^*\mu$ agree as forces along $L$
except at endpoints $x \in\bd L$, where they may differ by an
atomic force in the direction normal to $M_x$.
\qed \end{theorem}

Most of the configurations we care about here consist---as did the tight simple
clasp~[2]---of convex curves in perpendicular planes.
Let $P_i$ ($i=1,2$) be the $x_iz$-coordinate plane in $\R^3=\{(x_1,x_2,z)\}$.
If $\gamma_i\subset P_i$ are two components of our generalized link,
then we parametrize each curve $\gamma_i$ by
the $z$-component, $u_i$, of its unit tangent vector.
If there is a strut $\{p_1,p_2\}$ of length~$1$ connecting
these two components $\gamma_i$ then the curves are both
perpendicular to this strut, so elementary trigonometry
gives the following lemma from~[2]:

\begin{lemma}\label{lem:perp-planes}
Let $\gamma_1$ and~$\gamma_2$ be two components of a link $L$,
lying in perpendicular planes.  Suppose there is a strut $\{p_1,p_2\}$
of length~$1$ connecting these components.  Then in the notation
of the previous paragraph, the parameters $x_i$ and $u_i$ for the
points $p_i$ satisfy $0\le x_i\le u_i\le 1$,
and any two of the numbers $x_1,x_2,u_1,u_2$ determine the other two
(up to sign) according to the formulas
\begin{eqnarray*}
x_i^2 &=& 1-\frac{x_j^2}{u_j^2} = \frac{u_i^2(1-u_j^2)}{1-u_i^2u_j^2}, \\
u_i^2 &=& \frac{1-x_j^2/u_j^2}{1-x_j^2} = \frac{x_i^2}{1-x_j^2},
\end{eqnarray*}
where $j\neq i$.
\qed \end{lemma}

Note that the lemma above says nothing about balancing of forces,
but is merely a geometric fact about curves in perpendicular planes
that stay distance $1$ apart.
To balance symmetric planar curves, we make use of the following
lemma, again from~[2]:

\begin{lemma}\label{lem:vert-force}
Suppose $\gamma_i\in P_i$ is symmetric across the $z$-axis,
and parametrized by $u_i$ as above.
Consider the net curvature force of a mirror-image pair of
infinitesimal arcs of $\gamma$.  This net force
acts in the vertical direction with magnitude $2|\d u_i|$.
\qed \end{lemma}

\section{Weighted clasps}

\begin{figure}
\begin{center}\begin{overpic}{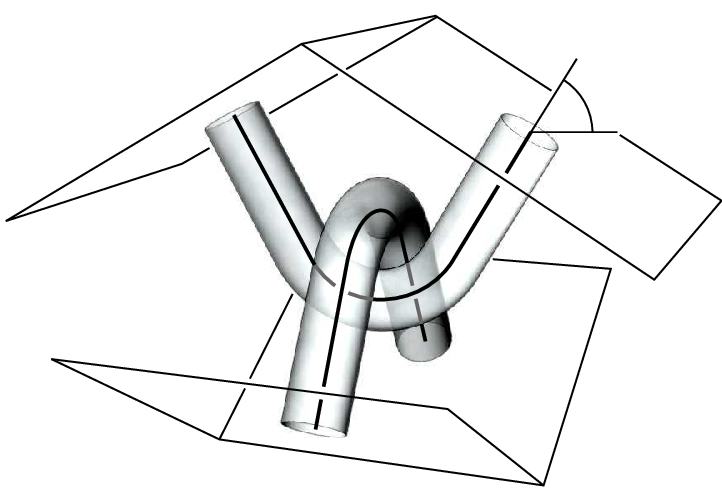}
\put(83,56){$\arcsin\tau$}
\end{overpic}\end{center}
\caption[The $\tau$-clasp problem]
{In the $\tau$-clasp problem,
the endpoints of two ropes are constrained to lie in four
planes whose normals make angle $\arcsin\tau$ with the horizontal.
The parameter $u_i$ ranges from $-\tau$ to $\tau$ along each arc,
as shown at the end of the top right arc. If extended, the four planes
shown would form the sides of a tetrahedron;
both arcs are constrained to stay within this tetrahedron.}
\label{fig:angleclasp}
\end{figure}

The {\em $\tau$-clasp} $C(\tau)$ is a generalized link
consisting of two clasped
ropes arranged according to the following description.  We fix four
planes, in a tetrahedral pattern as shown in \figr{angleclasp},
each making angle $\arcsin \tau$ with the vertical.
The $\tau$-clasp consists of two unknotted arcs $\alpha_1$ and $\alpha_2$ whose
endpoints are constrained to the four planes; the complement of
the tetrahedron also serves as an obstacle for the link. 
The isotopy class of the link is specified so that
closing each arc within the planes of its endpoints would produce a Hopf link.

If the two components $\alpha_1$ and $\alpha_2$ of a clasp have
different tensions,
the configuration can only be balanced if the opening angles
also differ.  For $\tau_1,\tau_2\in(0,1]$, we now define the weighted
clasp $C(\tau_1,\tau_2)$.  It is just like $C(\tau)$ except
that its component $\alpha_i$ is attached to planes
at angle $\arcsin\tau_i$.
For this generalized link, we will describe a critical configuration
for the weighted Gehring ropelength problem, where the weights
are $w_i=1/\tau_i$ on the two components.
These tensions are chosen to ensure a net balance of vertical forces
$2w_1\tau_1=2=2w_2\tau_2$ at the ends of the clasp.

It follows from \lem{vert-force} that the clasped contact between
$\alpha_1$ and $\alpha_2$ is determined
by the balance equation
$u_1/\tau_1 + u_2/\tau_2 = 1$.
Plugging this into \lem{perp-planes} gives us equations for the shapes.
Namely, we get an explicit formula for the $x_i$-coordinate of~$\alpha_i$:
$$x_i(u_i)
=\sqrt{\frac{u_i^2(1-u_j^2)}{1-u_i^2u_j^2}}
=\sqrt{\frac{u_i^2\big(1-\tau_j^2(1-u_i/\tau_i)^2\big)}
             {1-u_i^2\tau_j^2(1-u_i/\tau_i)^2}},$$
where $j\ne i$.
As for the critical $\tau$-clasp~[2],
the $z$-coordinate is then
determined as a hyperelliptic integral, through the relation
$\d z/\d x_i = u_i/\sqrt{1-u_i^2}$ that defines $u_i$.
This integral, and the similar one for the total arclength of
the curve, can easily be computed numerically.

Note that when $\tau_1=\tau_2$, the curves arising here are exactly
the symmetric $\tau$-clasp curves of~[2].
But when $\tau_1\ne \tau_2$, the two touching curves have shapes different
from each other and from any symmetric clasps.

The calculations above, combined with \thm{gen-gehr-bal}, serve to prove:
\begin{theorem} \label{thm:wtd-clasp}
The configuration of $C(\tau_1,\tau_2)$ described above
is critical for weighted ropelength.
\qed
\end{theorem}

We expect these configurations are in fact the ropelength minimizers.

\section{Clasps with parallels}
We next describe a family of examples based on
the connect sum of the $\tau$-clasp with
the Hopf link.  We describe configurations which we show
are balanced and thus critical points, but they are not minimizers
and are presumably very unstable equilibria.
Although we use unweighted length as our objective functional here,
some of the shapes that appear are those of the
weighted clasp curves we described above.

We first consider a configuration $C_{1,1}^1(\tau)$ which can be defined
as the connect sum of the $\tau$-clasp $C(\tau)$ with a Hopf link.
Letting $\beta$ and~$\gamma$ denote the two components of $C(\tau)$,
we add a third component
$\alpha$,  which is linked to $\gamma$ but not to $\beta$.
We guess that in the ropelength minimizer for this link,
the components $\beta$ and~$\gamma$
retain the shapes they have in the minimizing clasp, and $\alpha$ is a
round circle around a point on a straight end of $\gamma$, far from $\beta$.
Here, however, we will be interested in describing a critical
configuration with more symmetry, where $\alpha$ lies
in the plane of $\beta$.  

\begin{definition}
\label{defn:clink}
The generalized link $C_{k,l}^m(\tau)$ is defined from $C_{1,1}^1(\tau)$
by replacing $\alpha$ by $k$ parallel copies, $\beta$ by $l$ parallel copies,
and~$\gamma$ by $m$ parallel copies.  We also adjust the angles of the bounding
planes: those containing the ends of the $l$ copies of $\beta$ now
lie at angle $\arcsin(\tau/l)$ from the vertical, and those with $\gamma$
at angle $\arcsin(\tau/m)$.
\end{definition}

\begin{figure}
\begin{center}\begin{overpic}[width=1.5in]{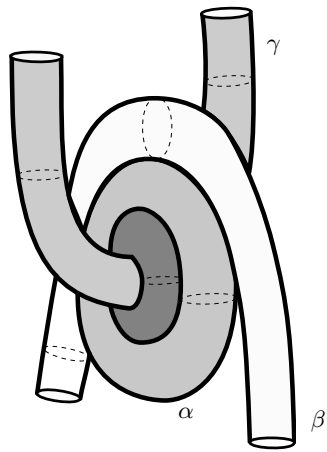}
\put(38,7){$\alpha$}
\put(68,5){$\beta$}
\put(58,90){$\gamma$}
\end{overpic}\end{center}
\caption[Clasp with linked loops]
{The generalized link $C_{2,1}^1(\tau)$ has a critical (but not minimizing)
configuration as shown here, with the two copies of $\alpha$ nested in
a common plane with~$\beta$.}
\label{fig:sw1}
\end{figure}

In the configurations we describe, the copies of $\alpha$ and $\beta$
all lie nested with one another in a single vertical plane $P_1$, and
the copies of $\gamma$ are nested in the perpendicular vertical plane $P_2$.
(\figr{sw1} sketches what such a configuration for
$C_{2,1}^1(\tau)$ might look like.)
In particular, the shapes of all components are determined by those of
the innermost $\alpha$ and~$\gamma$: the other copies of $\alpha$ are
successive outer parallels at distance~$1$ and the copies of $\beta$
are further outer parallels except that they peel off at the angle
corresponding to $u_1=\pm\tau/l$,
to proceed straight out and down to meet their bounding planes perpendicularly.
Similarly, the copies of $\gamma$ are determined as successive outer
parallels at distance~$1$ to the innermost copy; all of them have straight
segments out and up to the bounding planes at $u_2=\pm\tau/m$ and have a curved
arc parametrized by $u_2\in[-\tau/m,\tau/m]$.
Here $\tau$ can be any nonnegative number not exceeding $\min(m,n)$.
The angles are again determined by an overall balance of vertical forces:
the $l$ strands of $\beta$ exert a vertical force of
$2|u_1|=2\tau/l$ each, while the $m$ strands of $\gamma$ exert a 
vertical force of $2|u_2|=2\tau/m$ each.

Our entire configuration, like the $\tau$-clasp, has
mirror symmetry across the planes $P_1$ and $P_2$, but unless we have
$k=0$ and $l=m$, there is no longer the extra symmetry
interchanging top and bottom.

\begin{figure*}
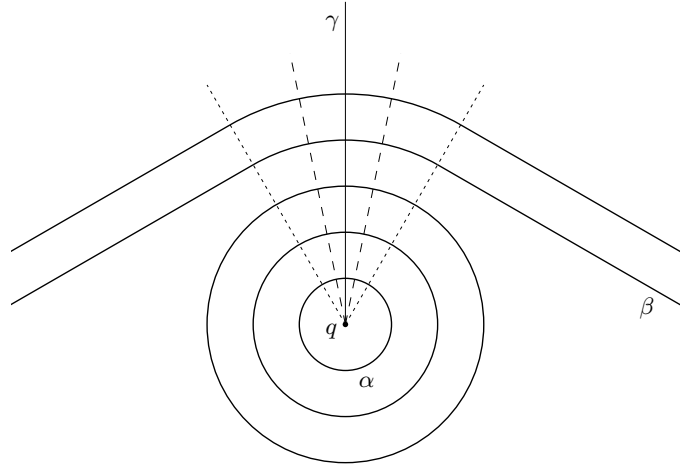

\begin{center}\begin{overpic}[width=3.5in]{mclasp}
\put(47,67){$\gamma$}
\put(47,21){$q$}
\put(52,13){$\alpha$}
\put(94,24){$\beta$}
\end{overpic}\end{center}
\caption[The multiple clasp problem]
{In this variant $C_{3,2}^1(1)$ of the clasp problem,
there are $k=3$ nested parallel copies
of a closed arc $\alpha$ inside $l=2$ parallel copies of a clasp arc $\beta$,
here proceeding out and down towards bounding planes below the figure.
These are linked to $m=1$ copies of a clasp arc $\gamma$ which
(in the case $\tau=1$ illustrated) is attached to a horizontal ceiling.
(We see $\gamma$ in projection to the plane $P_1$ of the figure only
as a vertical line; the dot at the center is its tip~$q$.)
The arcs $\beta$ peel away from the arcs $\alpha$ at $u_1=\tau/k$, at the
dotted line shown at $30^\circ$; the arcs $\alpha$ and $\beta$ are clasped to
$\gamma$ for $|u_1|\le\tau/(k+l)$ (above the dashed lines shown),
and then continue as circular arcs (between the dashed and dotted lines).
The figure has been drawn entirely with circles, but the true analytic
arcs near the tips would be imperceptibly different.
The arc~$\gamma$ is an analytic clasp for its entire curved
region $|u_2|\le 1$, and then (since $\tau/m=1$)
proceeds vertically to the ceiling.}
\label{fig:mclasp}
\end{figure*}

It remains to describe the shapes of the innermost copies $\alpha$
and~$\gamma$ exactly, and to show the resulting configuration is balanced.
(See \figr{mclasp}, which shows $C_{3,2}^1(1)$.)
We will find it useful to use the abbreviation $n:=k+l$ for the total
number of curves in the plane~$P_1$ containing $\alpha$ and~$\beta$
and their parallels.
Note that the $m$ copies of $\gamma$ intersect the plane $P_1$ in
a series of $m$ points $q_i$ spaced at distance $1$ along the $z$-axis.
The inner copy $\alpha$ would naively be a stadium curve looping
around these points at distance~$1$.  Our configuration is close
to this, changing only part of the upper semicircle.  That is, $\alpha$
consists of a semicircle around the bottommost $q_i$, joined to vertical
segments of length $m$, and to an upper arc parametrized by $u_1\in[-1,1]$.
For $|u_1|\ge\tau/n$ this upper arc is also part of the unit circle
around the topmost $q_i$, but for $u_1\in[-\tau/n,\tau/n]$ it
is an analytic arc determined below.  Note that the parallel copies
of $\beta$ include outer parallels of this analytic arc in the range
$|u_1|\le\tau/n$, but also outer parallels to the circle (that is,
larger circular arcs of integer radius around the topmost $q_i$)
in the range $|u_1|\in[\tau/n,\tau/l]$, before peeling off straight
at $|u_1|=\tau/l$.

The innermost $\gamma$ consists more simply of an analytic arc parametrized 
by $|u_2|\le\tau/m$ joined to the straight segments out to the bounding planes.
It now remains merely to describe the
analytic arcs of the innermost $\alpha$ and~$\gamma$: these form a weighted
clasp $C(\tau/n,\tau/m)$ as in \thm{wtd-clasp}.

The important observation is that in our situation we have a convex
planar curve with nested outer parallels.  Each of the parallels has struts
only to the next ones inward and outward, and these touch at points with
equal direction $u_i$.  Thus the innermost curve behaves like a curve
with increased tension, proportional to the total number of parallel curves.
If this increased force is balanced by struts to the other innermost curve,
then the struts between the parallels distribute this balancing force outwards
to balance the equal tension on each of the parallel strands.

\begin{theorem} \label{thm:single-simple-sum}
The configuration of $C_{k,l}^m(\tau)$ described above
is critical for ropelength.
\end{theorem}

\begin{proof}
The lower semicircle of~$\alpha$ and its parallels exert a total
force $2k$ upwards on the bottommost~$q_i$.  This is transmitted upwards
to the topmost~$q_i$ by atomic strut forces in the struts connecting
the~$q_i$.  The remaining circular pieces around this topmost~$q_i$,
namely $k$ circles for $|u_1|\in[\tau/l,1]$ and a total of~$n$ for
$|u_1|\in[\tau/n,\tau/l]$, exert a balancing downwards force of
$2k(1-\tau/l)+2n(\tau/l-\tau/n)=2k$ on this~$q_i$.
The analytic arcs of the innermost~$\alpha$ and~$\gamma$ stay
at constant distance~$1$ from each other; since they are the curves
of the critical weighted clasp $C(\tau/n,\tau/m)$, they balance
each other with weights~$n$ and~$m$, or equivalently, with unit weights
on each of them and their respective $n$ and~$m$ parallels.
\end{proof}

We now consider in more detail the specific example of $C_{1,1}^1(1)$.
This is the ordinary simple clasp~$C(1)$---whose endpoints are attached
to horizontal planes---with the addition of a closed component~$\alpha$.
The lower half of $\alpha$ is a semicircle of radius~$1$ centered at the
tip of $\gamma$. The upper half of $\alpha$ consists of three parts:
two circular arcs of angle $30^\circ$ and an arc clasped to
$\gamma$. The curve $\beta$ consists of two vertical segments from the
floor up to the height of the tip of $\gamma$, connected by an arc
that is an outer parallel to the upper half of $\alpha$.
The curve $\gamma$ includes no circular arcs, but only the analytic
arc determined above and straight segments up to the ceiling.

In general, note that our critical configuration of $C_{k,l}^m(\tau)$
includes circular arcs in the $\beta$ components
when $k>0$, but not when $k=0$.  For $k=0$ all the force balancing
happens between the analytic arcs described above.

\section{Conjectured minimizers for two cases}
\begin{figure*}
\centerline{
\begin{overpic}[height=2.3in]{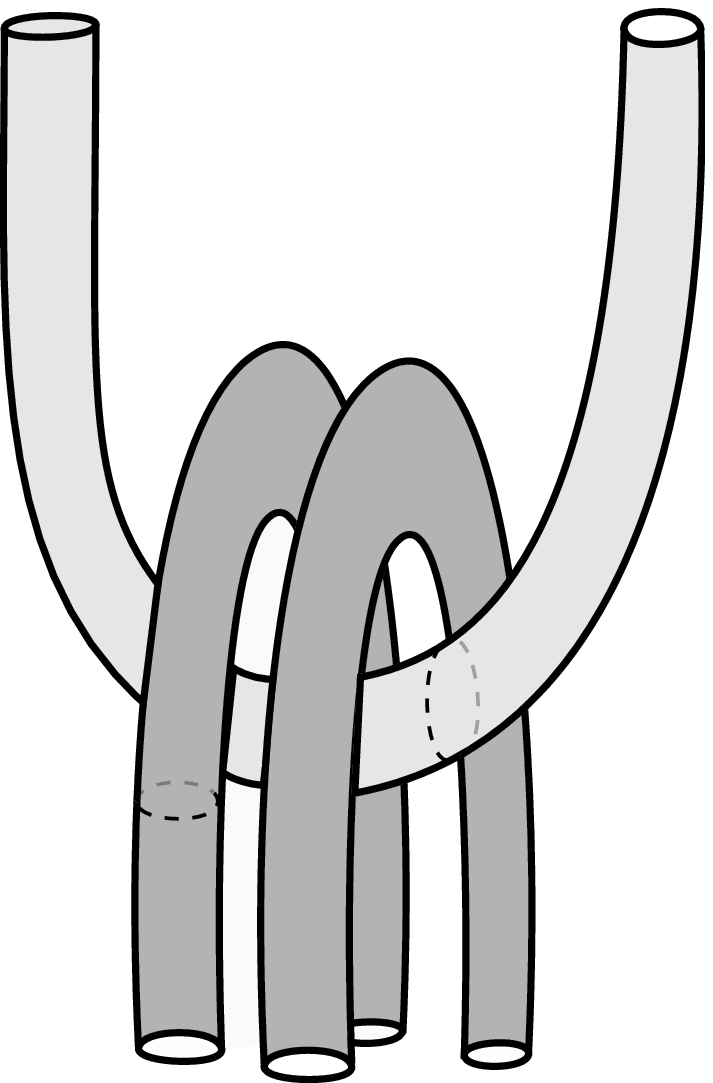}
   \put(53,5){$\beta$} \put(51,90){$\gamma$} \end{overpic}
\hfil
\begin{overpic}[height=2.3in]{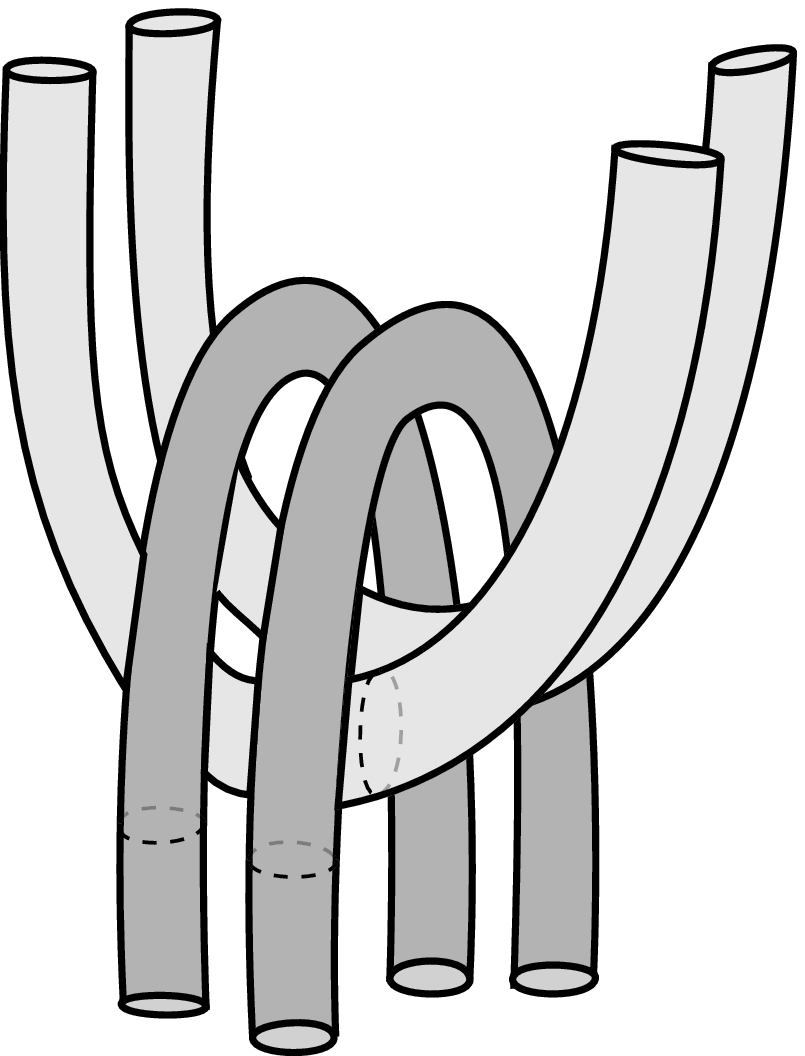}
   \put(60,10){$\beta$} \put(52,81){$\gamma$} \end{overpic}
}
\caption[The presumed minimzers for two generalized clasp problems]
{The configurations $L^m(\tau)$ are presumably the minimizers
for the clasp problems $C_{0,2}^m(\tau)$ for $m=1$ (left) and $m=2$ (right).
They involve congruent curves in parallel planes.
The exact geometry of the critical configurations---including
unit-length horizontal segments at the tips of most components, and the
fact that the lower ends on the left must proceed out within $30^\circ$ of
the horizontal---is not shown in these sketches.}
\label{fig:weave}
\end{figure*}

For the two cases $k=0$, $l=2$, $m=1,2$, we now describe a different
critical configuration $L^m(\tau)$ for the generalized link $C_{0,2}^m(\tau)$.
In both cases, we conjecture that this configuration (unlike the one
described above) is the minimizer for Gehring ropelength.
The drawings in \figr{weave} indicate the relative positions of
the components in these configurations, without showing the
exact geometric features described below.

For $C_{0,2}^1(\tau)$ the configuration we described above consists of
noncongruent curves $\beta_0$ and~$\gamma_0$, together with an outer parallel
to $\beta_0$.  In our new configuration $L^1(\tau)$, both $\beta$ curves
are congruent to $\beta_0$, but translated out to
lie in parallel planes $x_2=\pm\tfrac12$.  The $\gamma$ curve is a copy
of $\gamma_0$, but split apart at the tip, with a unit-length straight segment
inserted from $x_2=-\tfrac12$ to $x_2=\tfrac12$.

For $C_{0,2}^2(\tau)$ the configuration we described above consists of
two ordinary $\tau$-clasp curves together with their outer parallels.
In our new configuration $L^2(\tau)$, all four curves are congruent,
and each looks like the ordinary $\tau$-clasp but split at the tip, with
a unit-length segment inserted.  The two copies of $\beta$ lie in the
planes $x_2=\pm\tfrac12$, and are translates of each other perpendicular
to these planes.  Similarly the two copies of $\gamma$ lie in the planes
$x_1=\pm\tfrac12$.  Note that for $\tau=1$, the curvature of each curve
is unbounded near the tip;
$L^2(\tau)$ gives an example of a Gehring-critical configuration
in which curvature approaches infinity but then immediately jumps to zero.

\begin{theorem}
The configurations $L^1(\tau)$ and $L^2(\tau)$ described above
are critical for Gehring ropelength.
\end{theorem}

\begin{proof}
Consider the case $m=2$.  In the usual $\tau$-clasp, there is $2$--$2$ strut
contact, with the four struts in any given set lying over the four quadrants
of the $x_1x_2$-plane.  Here, those struts are pulled apart from each other.
They carry exactly the same balancing forces as before, but the horizontal
components of those forces need to be balanced by further horizontal struts.
These connect corresponding points on the two copies of $\beta$,
and similarly for $\gamma$.
Each quadrilateral of struts has now been replaced by an octagon
(with four translated copies of the original diagonal edges
plus four new horizontal edges).  But the forces still balance within each
such octagon.

The case $m=1$ is similar but even simpler.  In this case the quadrilateral
of struts is split apart only in one direction.  Again, new horizontal struts
connecting the copies of $\beta$ allow the horizontal components of force
to be properly balanced.  The struts here close into hexagons.
\end{proof}

\section{A chained clasp}
\begin{figure}
\centerline{\begin{overpic}[height=2.2in]{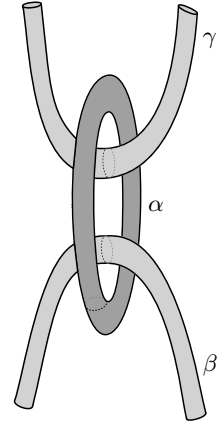}
 \put(32,50){$\alpha$} \put(45,12){$\beta$} \put(45,90){$\gamma$}
 \end{overpic}}
\caption[A chained clasp]
{Two clasps, when joined end to end, form a generalized link with three
components.  For the case $\tau=1$ of vertical ends,
the middle closed component~$\alpha$ has varying length
in our family of critical configurations.}
\label{fig:chain}
\end{figure}

Our next example is that of two clasps joined end to end,
as in \figr{chain}.
That is, we have three components: the open arcs $\beta$ and~$\gamma$ are
attached to the floor and ceiling, respectively, and are unlinked
with each other; each, however, is linked to the closed component~$\alpha$.
The configurations we consider have reflection symmetry across a
horizontal plane, interchanging $\beta$ and~$\gamma$ while preserving~$\alpha$.

If we make such a configuration with $\tau<1$, the junctions
(where~$\alpha$ clasps to~$\beta$ and to~$\gamma$) will move towards each
other, until the tips of~$\beta$ and~$\gamma$ touch at distance~$1$
producing an isolated strut.
Here $\beta$ and~$\gamma$ are congruent $\tau$-clasp curves.
The curve $\alpha$ is like a stadium curve, but its tips (for
$|u|\le\tau$) are $\tau$-clasps.  These are followed by circular
arcs around the tips of $\beta$ and~$\gamma$, which are finally
joined by unit-length vertical segments.
Assuming the top/bottom symmetry, this
critical configuration is uniquely determined.

On the other hand, when $\tau=1$, the critical configuration has
one simple clasp at the junction of $\alpha$ with~$\beta$, and
another symmetric one between $\alpha$ and~$\gamma$.  These clasps
can be close to each other (with the tips of $\beta$ and~$\gamma$
as close as distance $1$) or can move farther apart.  There is
a one-parameter family of equal-length symmetric critical configurations.
In this family, the length of the closed component $\alpha$ varies.

\begin{theorem}
The configurations of the chained clasps described above are critical
for Gehring ropelength.
\end{theorem}

\begin{proof}
For $\tau=1$ the balancing is just that for the $\alpha,\beta$
and $\alpha,\gamma$ clasps separately.

For $\tau<1$ we need to combine the balancing for the
clasp with that for the simple closed chain described in~[1,~2].
The curves~$\beta$ and~$\gamma$ are ordinary $\tau$-clasp curves,
positioned so their tips are unit distance apart.  These are balanced
by clasp arcs of~$\alpha$.  The circular arcs
of~$\alpha$ focus net force $2-2\tau$ downwards on the tip of~$\gamma$
and upwards on the tip of~$\beta$.  These forces are balanced by
a force on the isolated strut connecting these two tips.
\end{proof}

Note that we could build similar configurations with several chained
components~$\alpha_i$ in between~$\beta$ and~$\gamma$.  For $\tau<1$
the~$\alpha_i$ are all congruent to the curve~$\alpha$ described above.
For $\tau=1$, each would be a stadium curve, but they could have
differing lengths.  Again, we expect that all of these are ropelength
minimizers.

\section{The Granny Clasp}
\begin{figure}
\hfil \includegraphics[height=1.7in]{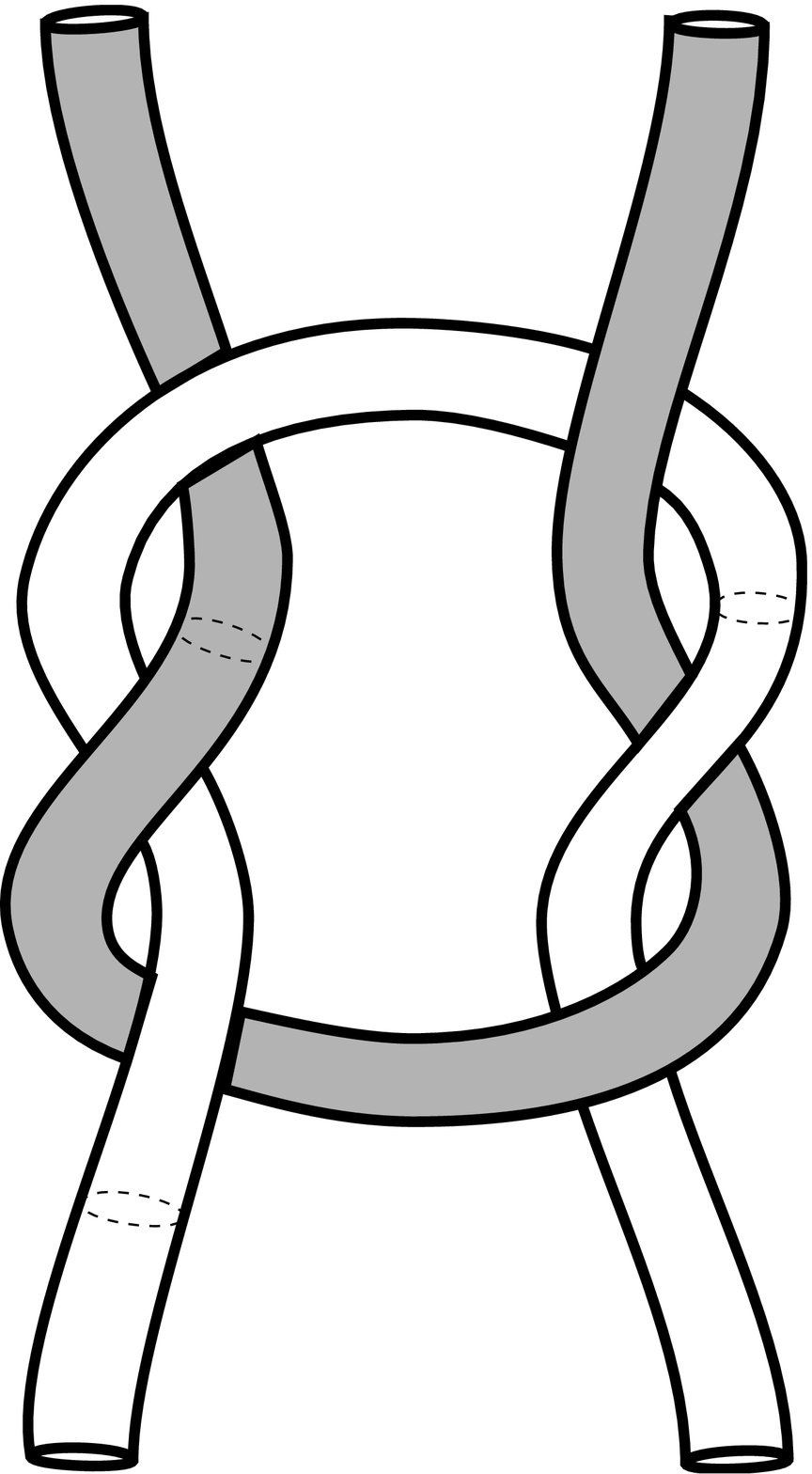}
\hfil \includegraphics[height=1.7in]{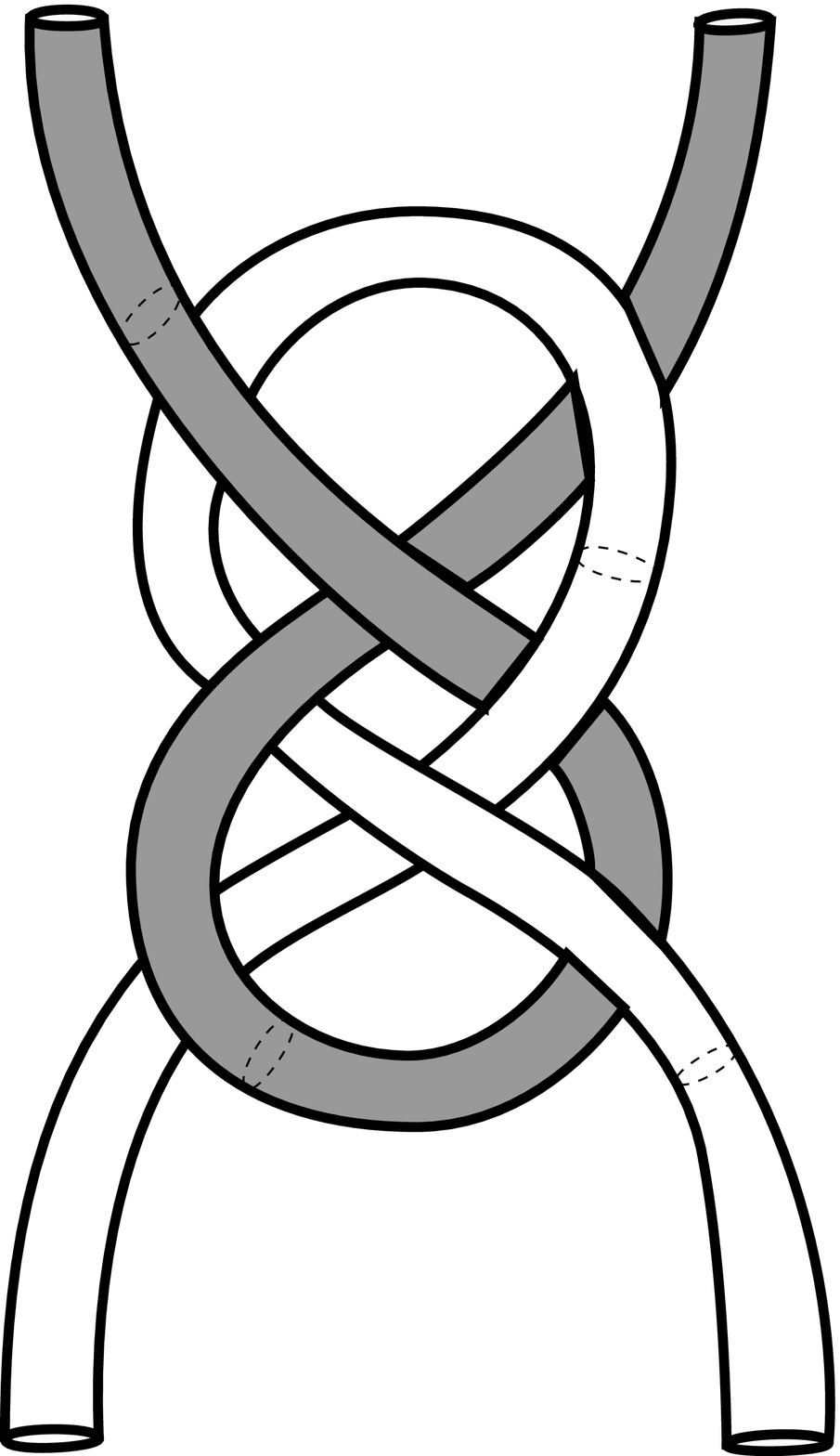}
\hfil \includegraphics[height=1.7in]{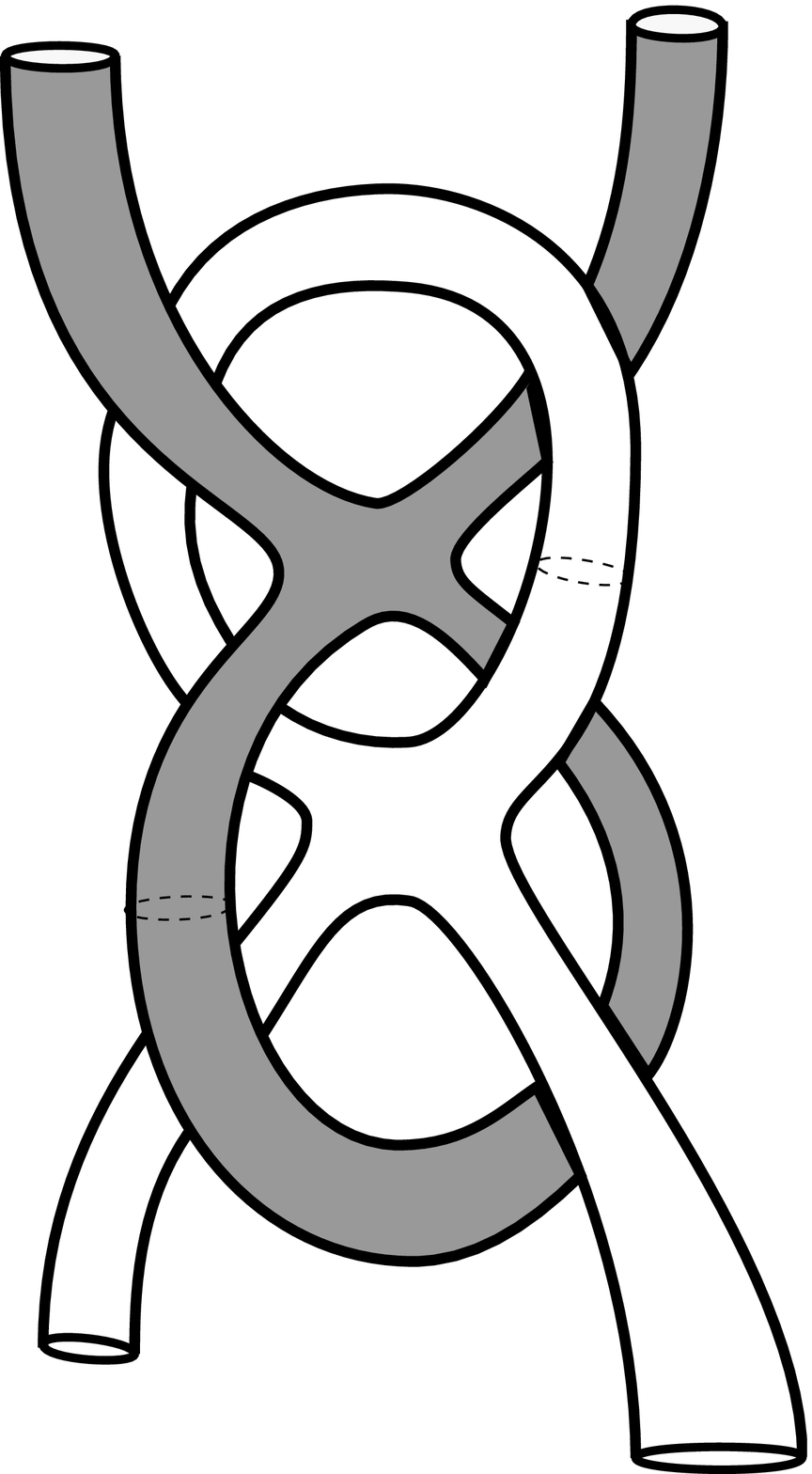}
\hfil
\caption[The granny clasp]
{This generalized link is called the granny clasp.
In the first configuration shown, if the two free ends at the
left were joined and likewise the two on the right, we would
get a granny knot, the connect sum of two trefoils of the same handedness.
In the middle, we see an isotopic configuration,
presumably close to a minimizer for the ordinary ropelength
problem~[1] where each component avoids itself as well as the others.
Note that, compared to the first picture, the upper two strands have
crossed, as have the lower two: the endpoints of the generalized link
are free to move within the floor and ceiling bounding planes.
For the Gehring ropelength we consider here, however,
the appropriate setting is Milnor's link-homotopy.
Our critical configuration (and presumed minimizer)
has components which are planar but fail to be
embedded curves, as on the right.}
\label{fig:granny}
\end{figure}

Our final example generalizs the simple clasp not by introducing extra
components, but by clasping the two components in a more intricate way.
The connect sum of two trefoil knots of the same handedness is
called the granny knot.
We define the granny clasp~$G_1$ to be
the generalized link shown in \figr{granny} (left),
a clasp of two ropes based on this granny knot.

It seems clear that if this configuration of the granny clasp were tied
tight in rope, each component would contact itself as well as the other,
as suggested in \figr{granny} (center).
But our constraint on the Gehring thickness $\gThi\ge1$ does not
see self-contact of a single component.  It is thus important to
remember that the natural setting for Gehring ropelength
problems~[2] is Milnor's link homotopy.
Two configurations are link-homotopic if there
is a homotopy between them where the components stay disjoint but
self-intersections of any given component are allowed.

\begin{figure}
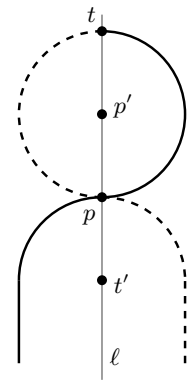

\begin{center}
\begin{overpic}[height=2in]{grannyhalf}
\put(27,25){$t'$}
\put(19,44){$p$}
\put(27,72){$p'$}
\put(20,96){$t$}
\put(26,6){$\ell$}
\end{overpic}
\end{center}
\caption[Critical granny]{In the critical configuration of the
granny clasp~$G_1$, each component is planar and built
from two symmetric halves, here shown
as solid and dashed lines.  Here $p$ is the point of self-intersection
and $t$ is the tip.  The other component is congruent to this one, and
lies in a perpendicular plane; the two places where it intersects
the plane of the figure are the points~$p'$ and~$t'$
corresponding to~$p$ and~$t$.  Each half of the component shown
consists of a straight arc from the floor up to the level of~$t'$,
followed by three $\tau=1$ clasp arcs, the first up to~$p$, the next to the
level of~$p'$ and the last up to~$t$.  (This figure has been drawn
with six quarter circles in place of the $\tau=1$ clasp arcs.)}
\label{fig:grannyhalf}
\end{figure}

In the critical configuration we describe for the granny clasp $G_1$,
each component does have a point of self-intersection, like those shown
in \figr{granny} (right).  More precisely, if~$\beta$ denotes the component
attached to the floor, it---like the curves in all our previous
examples---lies in a vertical plane and has mirror symmetry across
a vertical line $\ell$ in that plane.  In \figr{grannyhalf} we see
(as the solid line) one
symmetric half of~$\beta$, consisting of four analytic pieces
joined in a $C^1$ fashion: first a vertical segment up from the floor,
then a $\tau=1$ clasp arc leading to the point $p$ of self-intersection,
then continuing on the other side of $\ell$ with two more clasp arcs
leading to the tip $t$ of~$\beta$.  We expect that this configuration
is the minimizer, but as usual will prove only that it is critical.

\begin{theorem}
The (nonembedded) configuration of the generalized
link-homotopy class $G_1$ built from clasp arcs as described above
is critical for the Gehring ropelength problem.
\end{theorem}

\begin{proof}
The configuration can be balanced as follows: Cut space with horizontal
planes at the heights of~$t'$, $p$, $p'$ and~$t$.  Below~$t'$ and above
$t$ we have only vertical straight segments.  In each of the three intermediate
slabs we see exactly a $\tau=1$ clasp.  (Even though the symmetric halves
of~$\beta$ do not connect to each other through~$p$ as in a clasp, they still
balance the other component in the same way.)
\end{proof}

\begin{figure}
\hfil \includegraphics[height=2.1in]{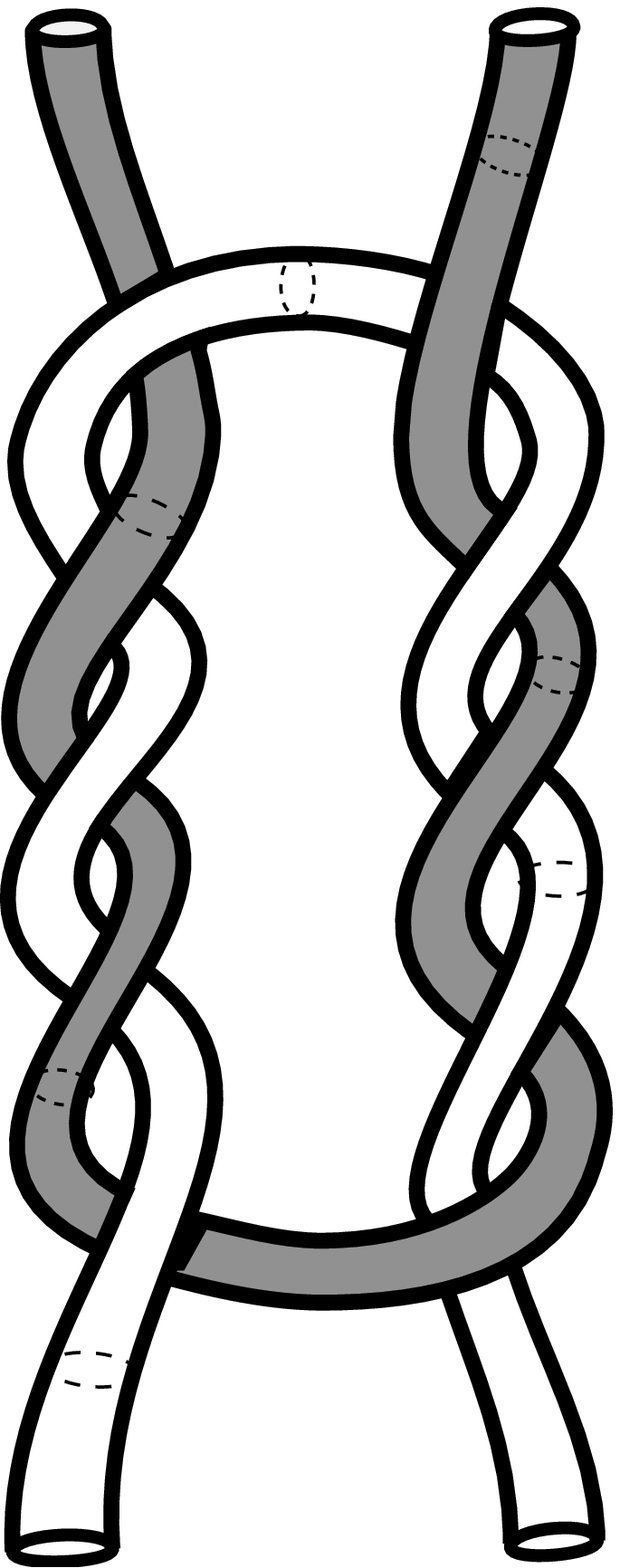}
\hfil \includegraphics[height=2.1in]{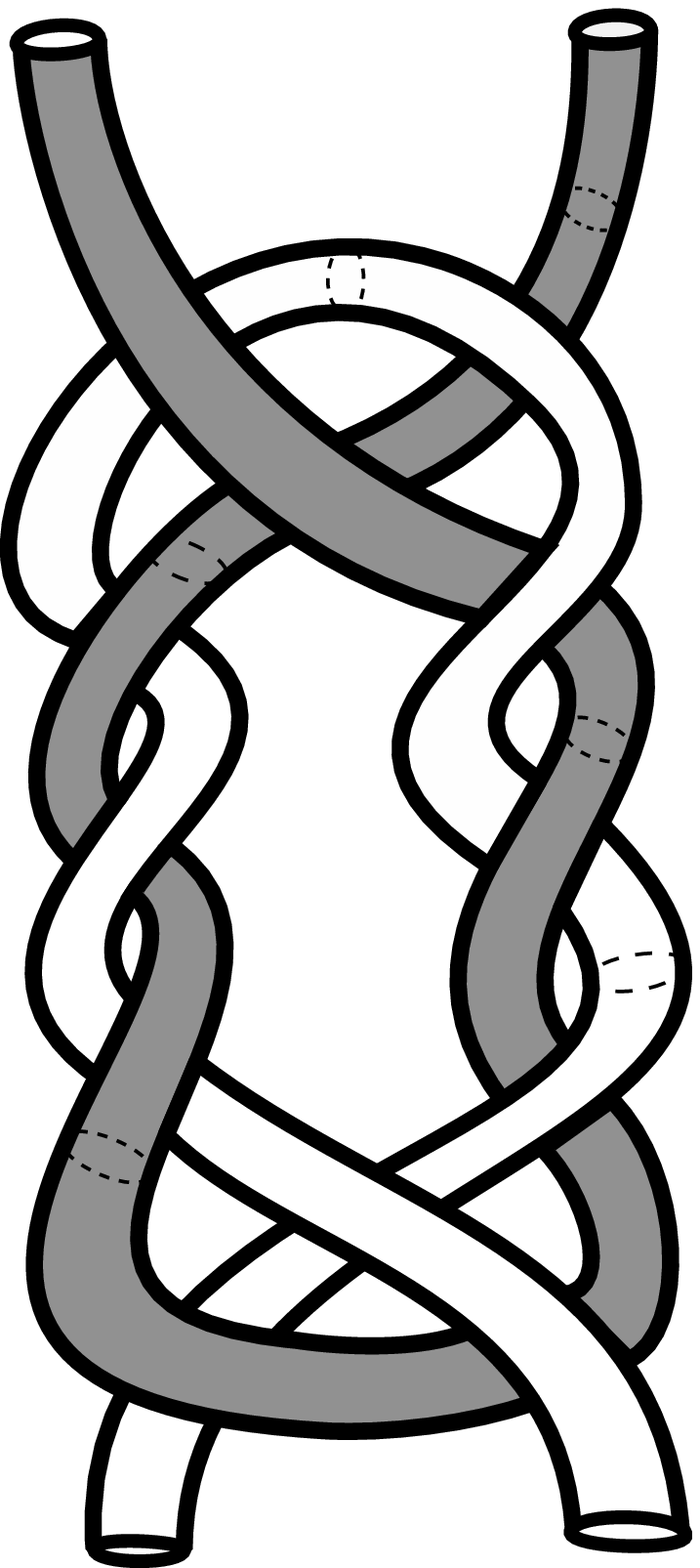}
\hfil \includegraphics[height=2.1in]{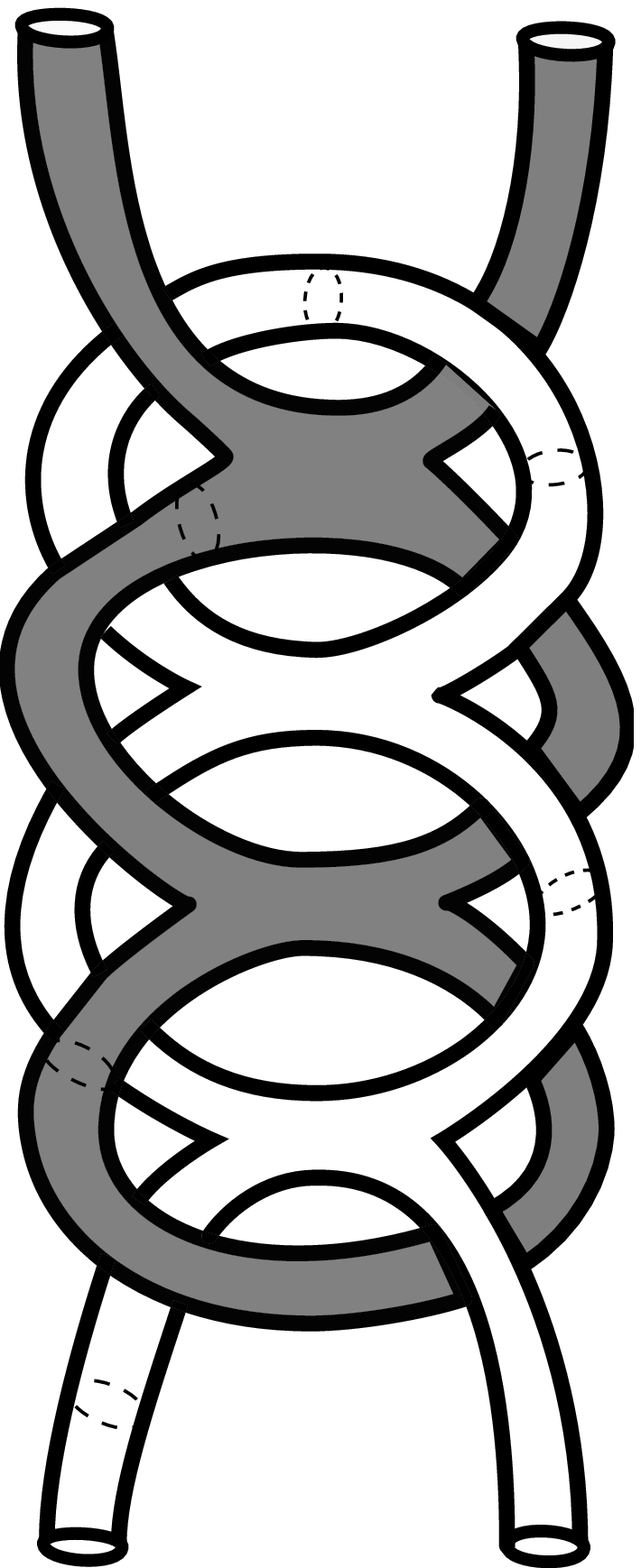}
\hfil
\caption[Higher-order granny clasp]
{This higher-order granny clasp $G_2$ is shown in three views.
It is derived~(left) from the connect sum of two $(2,5)$--torus knots.
Our conjectured
minimizer for the Gehring ropelength problem is a nonembedded
representative~(right) of this generalized link-homotopy class,
built from several clasp arcs.
The link homotopy from the initial configuration to this critical configuration
involves twisting the pairs of endpoints one full turn in opposite
directions.  In the center we see the configuration
after each has been twisted one half turn.
}
\label{fig:granny2} 
\end{figure}

We generalize this example as follows: let~$G_n$ be the generalized
link of two components obtained from the connect sum of two $(2,2n+1)$--torus
knots of the same handedness.  That is, each chain of three half-twists in
our first picture of $G_1$ is replaced by $2n+1$.  The link $G_2$ is
shown in \figr{granny2} (left); note that $G_0$ is the ordinary clasp $C(1)$.

Again we can describe a critical configuration for $G_n$---as shown
in \figr{granny2} (right)---which we expect is the mininimizer.
We obtain this configuration from the initial one by
twisting the pair of endpoints in the ceiling around each other $n$~full turns
(relative to the pair in the floor) and then letting
the $n$~points of self-contact of each component become self-intersections.
These $n$~self-intersections of each component occur
where it crosses its plane of symmetry. Each half of each component is built
from $2n+1$ clasp arcs plus a straight vertical segment.

A final generalization $G_n(\tau)$ would allow the ends of the clasp to be
attached to slanted planes.  Here the critical configuration would be
built from $\tau$-clasp arcs (for $|u|\le\tau$) near each self-intersection
point, and arcs of circles (for $|u|\ge\tau$) centered at the
self-intersections of the other component.

\section*{Acknowledgements}

We wish to thank Rob Kusner for helpful suggestions,
Jason Cantarella for the use of his $\tau$-clasp figure from~[2],
and TU~Berlin and the DFG Research Center {\textsc{Matheon}}
for hosting Wrinkle in January 2004 when this work was done.


\begin{thebibliography}{00}

\bibitem{CKS} J.~Cantarella, R.~Kusner, J.~Sullivan.
{\em On the minimum ropelength of knots and links}.
Inventiones Math. {\bf 150}:2, 2002, pp 257--286.
{\tt ArXiv:math.GT/0103224}.

\bibitem{CFKSW} J.~Cantarella, J.H.G.~Fu, R.~Kusner, J.~Sullivan, N.~Wrinkle.
{\em Criticality for the Gehring link problem}. Preprint, 2004.
{\tt ArXiv:math.DG/0402212}.

\bibitem{Sul} J.~Sullivan. {\em Curves of finite total curvature}.
Lecture notes from the May 2004 Oberwolfach Seminar
``Discrete Differential Geometry''.  Birkh{\"a}user, 2005, to appear.

\end{thebibliography}
\end{document}